\documentclass{article}
\usepackage[utf8]{inputenc}
\usepackage{amsmath}
\usepackage{graphicx}
\usepackage{physics}
\usepackage[ruled]{algorithm2e} 

\newcommand{\Obj}{O}
\newcommand{\x}{\mathbf{x}}
\renewcommand{\u}{\mathbf{u}}
\newcommand{\g}{\mathbf{g}}
\newcommand{\sfull}{\dv{\x}{\u}}
\newcommand{\s}{\mathbf{S}}
\newcommand{\Hess}{\mathbf{H}}

\date{}
\author{
  Simon Zimmermann
  \and
  Roi Poranne
  \and
  Stelian Coros
}

\title{Optimal control via second order sensitivity analysis}

\begin{document}

\maketitle
\noindent We show how to solve optimal control problems of the form
\begin{subequations}\label{eq:TheOptimizationProblem}
\begin{alignat}{2}
&\!\min_{u} &\qquad&  \label{eq:objective}
\Obj(\x,\u), \\
&\text{s.t.} &      & \g(\x,\u) = 0, \label{eq:gIsZero}
\end{alignat}
\end{subequations}
where $\x$ are \emph{state} variables, $\u$ are \emph{control} variables.
We use second order sensitivity analysis, and this short write-up presents a very sanitized derivation.
For a conceptually similar but more formal derivation, we refer the reader to~\cite{Jackson:1988:SecondOrderSensitivityAnalysis}.
The challenge of solving~\eqref{eq:TheOptimizationProblem} usually comes from the non-linearity of \eqref{eq:gIsZero}.
Ideally, \eqref{eq:gIsZero} would yield an explicit relationship $\x(\u)$, which would then be substituted into \eqref{eq:objective}, resulting in the \emph{un}constrained problem
\begin{equation}\label{eq:Unconstrained}
    \min_{u}  \qquad   \Obj(\x(\u),\u).
\end{equation}
However, this not generally the case, and in fact, not quite that necessary.
Indeed, the main reason for transforming \eqref{eq:TheOptimizationProblem} into \eqref{eq:Unconstrained} is for easily computing the derivatives of $\Obj(\x(\u),\u)$, to use with a gradient-based optimization.
However, the derivative can actually be readily computed using the implicit function theorem.

We begin by applying the chain rule on $\Obj(\x(\u),\u)$:
\begin{equation}	\label{eq:dOdp}
    \dv{\Obj}{\u} = \pdv{\Obj}{\x} \sfull + \pdv{\Obj}{\u},
\end{equation}
The term $\s := \sfull$ is known as the \emph{sensitivity} term.
The analytic expression for it can be found using the fact that $\g(\x,\u)$ is always zero (i.e. we assume that for {\em any} $\u$ we can compute $\x(\u)$ such that Eq.~\ref{eq:gIsZero} is satisfied), which implies:
\begin{equation}	\label{eq:dxdp}
    \dv{\g}{\u} = \pdv{\g}{\x} \s + \pdv{\g}{\u} = 0.
\end{equation}
By rearranging this equation, we get:
\begin{equation}\label{eq:SensitivityTerm}
    \s = - \left( \pdv{\g}{\x} \right)^{-1} \pdv{\g}{\u},
\end{equation}
and plugging into \eqref{eq:dOdp}, we obtain:
\begin{equation}	\label{eq:dOdp_complete}
    \dv{\Obj}{\u} = - \pdv{\Obj}{\x} \left( \pdv{\g}{\x} \right)^{-1}  \pdv{\g}{\u} + \pdv{\Obj}{\u}.
\end{equation}
We note that through a reordering of matrix multiplications, the well-known adjoint method avoids computing $\sfull$ directly as it evaluates $\dv{\Obj}{\u}$. This is oftentimes more computationally efficient.
However, we can leverage $\sfull$ to derive a second-order solver that exhibits much better convergence properties than first order alternatives.

To begin with, we differentiate \eqref{eq:dOdp}:
\begin{equation}\label{eq:tmpd2dp2}
    \dv[2]{\Obj}{\u} = \dv{\u} \dv{\Obj}{\u} =
    \dv{\u}( \pdv{\Obj}{\x}\s ) +
    \dv{\u} \pdv{\Obj}{\u}.
\end{equation}
The formulas above involve third-order tensors, which lead to notation that is slightly cumbersome. For conciseness, we treat tensors as matrices and assume that contractions are clear from context. 
The second term in Eq.~\ref{eq:tmpd2dp2} is straightforward:
\begin{equation}
    \dv{\u} \pdv{\Obj}{\u} = \s^T \pdv{\Obj}{\x}{\u}  + \pdv[2]{\Obj}{\u},
\end{equation}
while the first term evaluates to
\begin{equation}
    \dv{\u}( \pdv{\Obj}{\x}\s ) = \left(\dv{\u} \pdv{\Obj}{\x}\right)\s + \pdv{\Obj}{\x} \left(\dv{\u}\s\right),
\end{equation}
with
\begin{equation}
    \dv{\u} \pdv{\Obj}{\x} = \s^T \pdv[2]{\Obj}{\x} + \pdv{\Obj}{\x}{\u}.
\end{equation}
Here, $\dv{\u}\s$ is a third-order tensor, and $\pdv{\Obj}{\x} \left( \dv{\u}\s \right)$ stands for
$$
\pdv{\Obj}{\x} \left( \dv{\u}\s \right) = \sum_i \pdv{\Obj}{x_i} \left(\dv[2]{x_i}{\u}\right).
$$
The second-order sensitivity term $\dv{\u} \s$ must be further broken down:
\begin{equation}
    \dv{\u} \s  = \left( \s^T \pdv{\x}\s + \pdv{\u} \s \right).
\end{equation}
The partial derivatives of $\s$ can be found by taking the second derivatives in \eqref{eq:dxdp} and rearranging the terms.
This results in
\begin{equation}
    \pdv{\x} \s =
    - \left( \pdv{\g}{\x} \right)^{-1} \left(
    \pdv[2]{\g}{\x}\s
    + \pdv{\g}{\u}{\x} \right),
\end{equation}
\begin{equation}
    \pdv{\u} \s =
    - \left( \pdv{\g}{\x} \right)^{-1} \left(
    \pdv{\g}{\x}{\u}\s
    + \pdv[2]{\g}{\u} \right),
\end{equation}
where once again we assume that the tensor expressions are self-evident.
Combining all of the terms above leads to the following formula for the Hessian:
\begin{equation}\label{eq:TrueHessian}
    \dv[2]{\Obj}{\u} = 
    \pdv{\Obj}{\x} \left( \s^T \pdv{\x}\s + \pdv{\u}\s \right) + 
    \s^T \left( \pdv[2]{\Obj}{\x} \s +
    2 \frac{\partial^2 \Obj}{\partial \x \partial \u} \right) +  
    \frac{\partial^2 \Obj}{\partial \u^2}.
\end{equation}

\paragraph{Generalized Gauss-Newton}
Although Newton's method generally converges much faster than L-BFGS or gradient descent, there are two issues with it.
First, evaluating the second-order sensitivity term takes a non-negligible amount of time.
Second, the Hessian is often indefinite, and needs to be regularized.
Both problems can be dealt with by simply excluding the tensor terms in \eqref{eq:TrueHessian}.
The result is a generalized Gauss-Newton approximation for the Hessian:
\begin{equation}\label{eq:GaussNewton_d2OdP2}
    \Hess = 
    \s^T \pdv[2]{\Obj}{\x} \s +
    2 \s^T \pdv{\Obj}{\x}{\u}   +  
    \pdv[2]{\Obj}{\u}.
\end{equation}
Although $\Hess$ is not guaranteed to be positive-definite, we note that in many cases, $\Obj$ is a convex function of $\x$ and $\u$, and therefore the first and last terms are always positive definite. Additionally, $\Obj$ commonly does not explicitly couple $\x$ and $\u$, and therefore the mixed derivative (i.e. the second) term vanishes, which means that overall $\Hess$ is positive definite.

\paragraph{Sensitivity analysis for trajectories}
The implicit relationship described by Eq.~\eqref{eq:gIsZero} is very general, and can easily be derived for different types of dynamical systems. Turning directly to the time-discretized setting, the state vector $\x$ has dimension $n T$, where $T$ is the number of time steps and $n$ is the number of variables representing the state at one time step. In this setting, we let $\x_i$, the $i$-th $n$-block in $\x$, denote the configuration of the system at time $t_i$. The dynamical system obeys the \emph{state equation}
\begin{equation}\label{eq:equationsOfMotion}
\g\left(\x_i,\dot{\x}_i,\ddot{\x}_i,\u_i\right) = 0
\end{equation}
where $\dot{\mathbf{x}}_i,\ddot{\mathbf{x}}_i$ are time-discretized velocity and acceleration. For a simple, first order approximation, $\dot{\x}=(\frac{\x_i-\x_{i-1}}{h})$ and $\ddot{\x}=(\frac{\x_i-2\x_{i-1}+\x_{i-2}}{h^2})$, and so \eqref{eq:equationsOfMotion} can also be written as
\begin{equation}
\g\left(\x_i,\x_{i-1},\x_{i-2},\u_i\right) = 0
\end{equation}
Starting from an input {\em control} trajectory $\u$, and two fixed configurations, $\x_{0}$ and $\x_{-1}$, the entirety of $\x$ can be via forward simulation.

\begin{figure}
	\centering
    \includegraphics[width=0.7\linewidth]{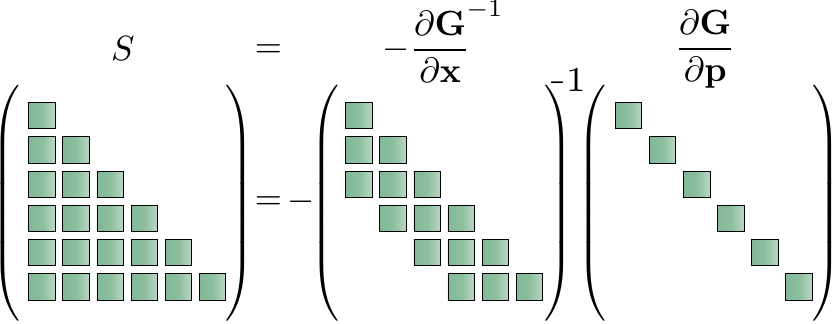}
	\caption{The structure of the system can be used to compute $\s$ faster by block forward-substitution.}
	\label{fig:blockSolver}
\end{figure}

The standard numerical integration scheme also reveals the structure of the implicit relationship described by Eq.~\eqref{eq:gIsZero}.
It is easy to see that the time domain imposes a very specific structure on the system of equations that must be solved to compute the Jacobian $\sfull$ in Eq.~\ref{eq:dxdp}. This structure is visualized in Fig.~\ref{fig:blockSolver}, and can be exploited to speed up computations.
In particular, since $\g$ depends explicitly only on $\x_i$, $\x_{i-1}$, $\x_{i-2}$ and $\u_i$ (i.e. all other partial derivatives are 0), $\pdv{\g}{\u}$ has a block diagonal form, and $\pdv{\g}{\x}$ has a \emph{banded} block diagonal form.
This allows us to solve the resulting system using block forward-substitution, rather than storing and solving the entire linear system represented by $\pdv{\g}{\x}$. We also note that the resulting $\s$ is block triangular, which correctly indicates that $\x_i$ does not depend on $\u_j$ if $j>i$, or intuitively, the control at any moment in time only affect future states.

\paragraph{Iterative optimization}
Using either $\dv[2]{\Obj}{\u}$ or its approximation $\Hess$, we can minimize \eqref{eq:TheOptimizationProblem} using a standard unconstrained optimization scheme, but we note one key difference: for each candidate $\u$, we must always compute the corresponding $\x$ to ensure that \eqref{eq:gIsZero} holds before evaluating any of the derivatives. Once $\pdv{\Obj}{\u}$ and $\Hess$ (or $\dv[2]{\Obj}{\u}$) are computed, the search direction $\mathbf{d}$ is found by solving:
\begin{equation}
    \mathbf{\Hess} \mathbf{d} = - \dv{\Obj}{\u}.
    \label{eq:searchDir}
\end{equation}
We use a backtracking line search to find the step size $\alpha$, where again, $\x$ need to be recomputed for each test candidate $\overline{\u}=\u+\alpha\mathbf{d}$ in order to evaluate $\Obj(\x(\overline{\u}),\overline{\u})$.
We summarize the optimization procedure in Algorithm \ref{alg}.

\begin{algorithm}[h!] 
  \caption{Trajectory optimization}
  \KwIn{Dynamical system, initial $\u$, initial $\x_0,\dot{\x}_0,$}
  \KwOut{Optimal control trajectory $\u$}
  \phantom{A} 
\While{criterion not reached} {
Compute $\x(\u)$ using forward simulation\\
Compute $\dv{\Obj}{\u}$ (Eq. \eqref{eq:dOdp}) \\
Compute $\mathbf{\Hess}$ (\eqref{eq:TrueHessian} or  \eqref{eq:GaussNewton_d2OdP2})\\
Solve $\mathbf{\Hess} \mathbf{d} = - \dv{\Obj}{\u}$\\
Run backtracking line search in $d$\\
\tcc{Simulate after every line search iteration}
} \label{alg}
\end{algorithm}

\end{document}